\documentclass[journal]{IEEEtran}

\ifCLASSINFOpdf
\else
   \usepackage[dvips]{graphicx}
\fi
\usepackage{url}

\usepackage{graphicx}
\usepackage{amssymb}
\usepackage{amsmath}
\usepackage{amsthm}
\usepackage{hyperref}
\usepackage{tabularx}
\usepackage{array}
\usepackage{booktabs}
\usepackage{subfig}
\usepackage{algorithm}
\usepackage{algpseudocode}
\usepackage{orcidlink}

\usepackage[utf8]{inputenc}
\DeclareUnicodeCharacter{03B4}{\ensuremath{\delta}}
\usepackage[shortlabels]{enumitem}
\newtheorem{thm}{Theorem}
\newtheorem{lem}[thm]{Lemma}   
\theoremstyle{definition}
\newtheorem{definition}{Definition}

\begin{document}

\title{Spectral Contraction of Boundary-Weighted Filters on $\delta$-Hyperbolic Graphs}

\author{%
  Le Vu Anh\,\orcidlink{0009-0000-1904-5186},\thanks{Le Vu Anh is with the Institute of Information Technology, Vietnam Academy of Science and Technology, 18 Hoang Quoc Viet, Ha Noi 100000, Vietnam (e-mail: levuanh@ioit.ac.vn).}%
  \and
  Mehmet Dik\,\orcidlink{0000-0003-0643-2771},\thanks{Mehmet Dik is with the Department of Mathematics, Computer Science \& Physics, Rockford University, 5050 E State St, Rockford, IL 61108, USA (e-mail: mdik@rockford.edu).}%
  \and
  Nguyen Viet Anh\,\orcidlink{0000-0001-7736-2470}\thanks{Nguyen Viet Anh is also with the Institute of Information Technology, Vietnam Academy of Science and Technology, Ha Noi, Vietnam (e-mail: anhnv@ioit.ac.vn).}%
}

\maketitle

\begin{abstract}
Hierarchical graphs often exhibit tree-like branching patterns, a structural property that challenges the design of traditional graph filters. We introduce a boundary-weighted operator that re-scales each edge according to how far its endpoints drift toward the graph’s Gromov boundary. Using Busemann functions on $\delta$-hyperbolic networks, we prove a closed-form upper bound on the operator’s spectral norm: every signal loses a curvature-controlled fraction of its energy at each pass. The result delivers a parameter-free, lightweight filter whose stability follows directly from geometric first principles, offering a new analytic tool for graph signal processing on data with dense or hidden hierarchical structure.
\end{abstract}


\IEEEpeerreviewmaketitle

\section{Introduction}\label{sec:intro}

The past decade has seen a growing interest in graph signal processing (GSP), where the goal is to analyze data supported on irregular domains, such as networks and graphs \cite{li2023gspreview}. Deep variants, such as Graph Neural Networks (GNNs), now appear in every branch of computational science and engineering, from materials discovery to epidemic forecasting \cite{Reiser2022GNNMaterials, Hy2022TemporalGNN}. This implementation has exposed key geometric mismatches between standard GNN architectures and the hierarchical nature of many real-world graphs. It motivates a search for alternative approaches based on non-Euclidean geometry. Among these, hyperbolic geometry has proven to be an effective tool for modeling graphs with dense hierarchical or tree-like structure \cite{zhu2024eventcentric}.

Many real-world graphs demonstrate exponential expansion, such as taxonomies, call graphs, and transportation hierarchies. These graphs tend to be nearly tree-like in structure, making hyperbolic geometry a natural choice for modeling. Unlike Euclidean embeddings, hyperbolic representations can encode such hierarchies with significantly lower distortion. As a result, hyperbolic GNNs have demonstrated superior accuracy, robustness, and sample efficiency. This advantage has been validated in recent studies on few-shot learning \cite{choudhary2023hgram}, temporal link prediction \cite{bai2023hgwn}, and deep architectures designed to mitigate the over-smoothing effect \cite{liu2024deephgcn}. The combination of spectral filtering with hyperbolic embeddings, such as via wavelets and multi-scale designs, also strengthens this trend \cite{yu2025hybowavenet}.

Beyond purely neural methods, a new class of spectral and physics-inspired models has been considered. One recent example is the work of Liu et al. (2024), which introduces a continuous geometry-aware graph diffusion framework based on hyperbolic neural partial differential equations \cite{Liu2024_ECML}. By formulating message passing as a continuous-time diffusion process in hyperbolic space, the model leverages the natural geometric alignment of hierarchical graph data. In addition to learned representation, stability has been a priority for the recent GSP technique design, such as how the blind deconvolution mechanism aligns with graph perturbations \cite{ye2025blindSPL}. In parallel, adaptive loss functions from signal filtering have been successfully transferred to hyperbolic GSP tasks \cite{cai2023lhcaf}.

Despite these developments, nearly all existing message-passing mechanisms remain restricted by operations within the hyperbolic manifold itself. They neglect directional cues from the Gromov boundary, a structure that governs signal behavior at large geodesic distances. So far, these asymptotic properties have only appeared indirectly, like in differentially private mechanisms \cite{wei2024poindp}, deep residual architectures \cite{xue2024rhgcn}, or open-world category discovery \cite{liu2025hypcd}, but not as clear edge-level weights in a GSP framework. In short, the following question remains open: \textit{How can boundary directions be explicitly embedded into local message-passing operations?}

We address this problem by proposing a boundary-weighted graph filter that scales each edge via a single exponential Busemann factor. This construction directly adds Gromov-boundary information without increasing the model’s parameter count. The resulting linear operator satisfies a spectral-norm bound of the form $\lVert T \rVert_2 \le e^{-\alpha\delta}$ for any $\delta$-hyperbolic graph, where $\alpha$ is a curvature-dependent scaling parameter. 

In summary, our key contributions are as follows:
\begin{itemize}
    \item \textbf{Efficient boundary-aware edge weighting.} We derive a vectorised, parameter-free Busemann weighting that augments graph filters with asymptotic geometry without increasing complexity.
    \item \textbf{Curvature-dependent stability guarantee.} We prove that the operator contracts the $\ell^2$-energy of any input signal by at least $e^{-\alpha \delta}$, establishing \textbf{the first explicit spectral bound} linked to the graph's hyperbolicity constant.
\end{itemize}

In the remainder of this letter, Section~\ref{sec:prelim} formalizes the graph-theoretic and geometric setup. Section~\ref{sec:operator} introduces the proposed boundary-weighted operator and establishes its contraction property as a function of graph curvature. Section~\ref{sec:concl} concludes with discussions and directions for further theoretical development.

\section{Problem setup}\label{sec:prelim}
We start by formalizing the negatively curved setting that initiates the proposed boundary-weighted filter.

\begin{figure}[h]
  \centering
  \includegraphics[width=\linewidth]{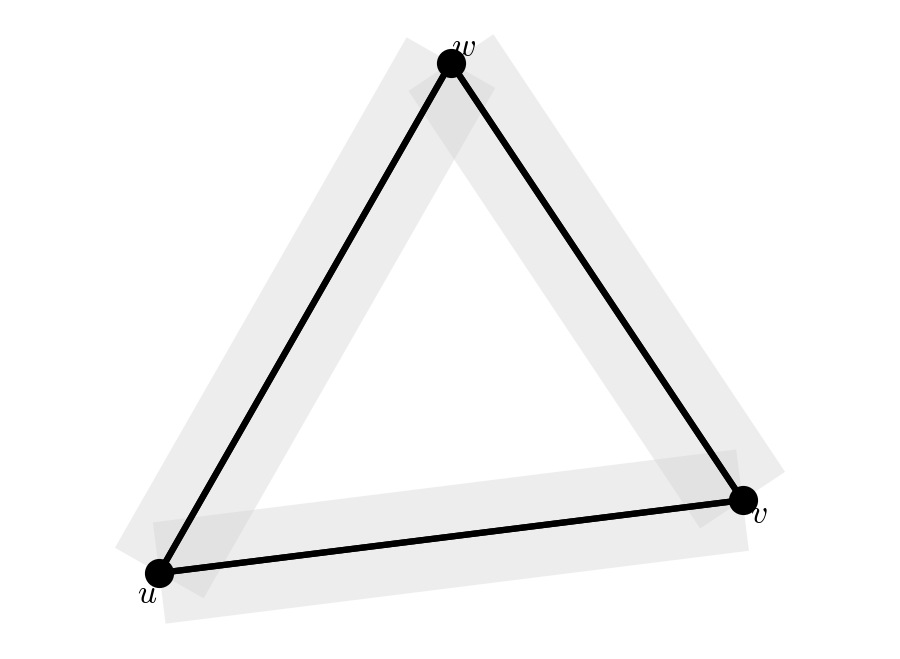}%
  \caption{A geodesic triangle in a $\delta$-hyperbolic graph. 
           The shaded corridor of width~$\delta$ shows that each side stays within distance $\delta$ of the union of the other two. This shaded region shows the tree-like geometry that initiates our
           boundary-weighted filter.}
  \label{fig:delta_triangle}
\end{figure}

\begin{itemize}[leftmargin=1.8em,itemsep=3pt]

  \item \textbf{Graph $G$.} We work with a finite, connected, undirected graph $G=(V,E)$ and write $d_G(u,v)$ for the usual shortest–path distance.

  \item \textbf{$\delta$-hyperbolicity.}  The graph is $\delta$-hyperbolic if every geodesic triangle is $\delta$-thin. It means each side lies inside the $\delta$-neighbourhood of the two others. Figure~\ref{fig:delta_triangle} shows such a triangle with vertices $u,v,w$. The shaded corridor of width~$\delta$ visualizes the tree-like "thinness" that provides a base for our contraction bound.

  \item \textbf{Gromov boundary $\partial G$.}  Points in $\partial G$ are equivalence classes of geodesic rays that stay a bounded Hausdorff distance apart. They represent directions to infinity and embed the asymptotic geometry of~$G$.

  \item \textbf{Anchor pair $(o,\xi)$.}  We fix a base vertex $o\!\in\!V$ and a boundary point $\xi\!\in\!\partial G$.  All boundary-aware quantities are measured relative to this reference pair.

  \item \textbf{Busemann function $\beta_{o,\xi}$.}  For every vertex $v\!\in\!V$
    \[
    \beta_{o,\xi}(v)=
    \lim_{t\to\infty}\Bigl[d_G\bigl(v,\gamma_\xi(t)\bigr)-t\Bigr],
    \]
    where $\gamma_\xi$ is any geodesic ray from $o$ toward~$\xi$. The limit exists and is independent of the chosen ray. As a result, a 1-Lipschitz “height" of $v$ is given along the boundary direction.

  \item \textbf{Scale parameter $\alpha>0$.}  A positive scalar controls how strongly boundary discrepancies attenuate edge messages. In section~\ref{sec:operator} we show that $\alpha$ also manages the operator’s spectral contraction.

\end{itemize}

\section{Boundary-weighted filter}\label{sec:operator}

\subsection{Main construction}\label{sec:bwf}

We now translate the geometric basis from section~\ref{sec:prelim} into a drop-in edge reweighting rule that can be inserted into any message-passing layer.

\begin{definition}[Boundary-weighted edge factor]\label{def:edgefactor}
For every edge $(u,v)\in E$, define
\begin{equation}\label{eq:edgeweight}
  w_{\alpha}(u,v)\;=\;
  \exp\!\Bigl[-\alpha\,\bigl|\beta_{o,\xi}(u)-\beta_{o,\xi}(v)\bigr|\Bigr],
\end{equation}
where $\alpha>0$ tunes how aggressively boundary divergence is penalized, and $(o,\xi)$ is a fixed anchor pair as defined in Section~\ref{sec:prelim}.
\end{definition}

\begin{figure}[h]
  \centering
  \includegraphics[width=\linewidth]{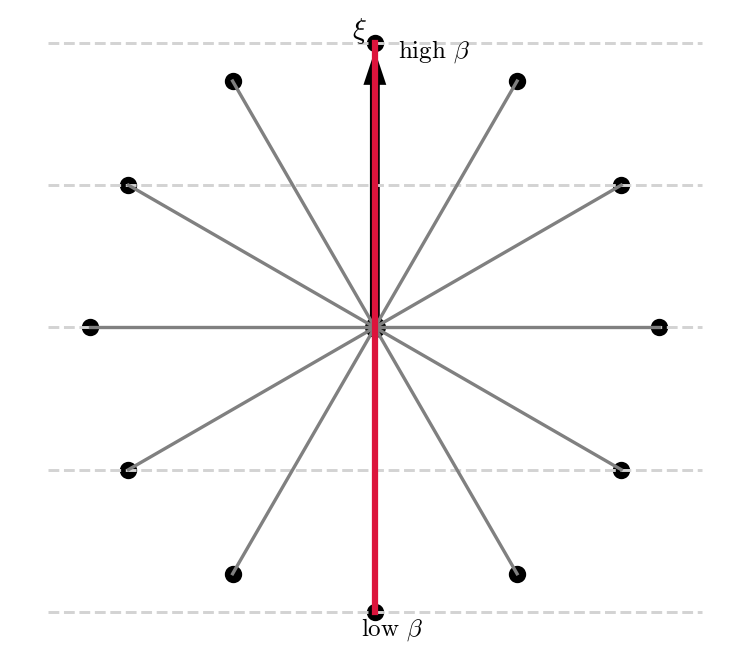}%
  \caption{Busemann "height" in a $\delta$-hyperbolic graph.  
           Concentric contour lines mark equal values of $\beta_{o,\xi}$, measured from the base vertex $o$ toward the boundary ray $\xi$ (arrow). The highlighted edge connects a low-$\beta$ node to a high-$\beta$ node and is therefore heavily down-weighted by the boundary factor in definition \ref{def:edgefactor}.}
  \label{fig:busemann_contours}
\end{figure}

Since $\beta_{o,\xi}$ increases at unit speed along every geodesic ray toward $\xi$, the quantity $|\beta_{o,\xi}(u)-\beta_{o,\xi}(v)|$ can be read as the vertical height
gap between the two endpoints of an edge when the graph is unfolded along the boundary ray. Definition~\ref{def:edgefactor} therefore assigns large weights to edges that run \emph{horizontally} within the same
contour band and small weights to those that jump across bands.  In fig.~\ref{fig:busemann_contours}, the dashed horizontal lines mark equal height levels. The highlighted edge crosses several levels, so its weight is strongly attenuated.  Conversely, edges lying entirely inside one band
retain almost their full influence, which is necessary for preserving the low-frequency structure encoded by the hierarchy.

\begin{definition}[Boundary-weighted layer]\label{def:bwlayer}
Given node features $f:V\!\to\!\mathbb{R}^{d}$ and a unit-norm mixing matrix $A\in\mathbb{R}^{d\times d}$, set
\begin{equation}\label{eq:Talpha}
  (T_{\alpha}f)(v)\;=\!\!\!
  \sum_{u\in N(v)} w_{\alpha}(u,v)\,A\,f(u),
  \qquad v\in V,
\end{equation}
\end{definition}

\smallskip
\textbf{Computational cost.}
We assume that the vertex–wise Busemann values $\{\beta_{o,\xi}(v)\}_{v\in V}$ are pre-computed once in $O(|V|)$ time.  Evaluating~\eqref{eq:Talpha} for a single edge $(u,v)$ includes 
\begin{itemize}
    \item (i) one subtraction and absolute value,
    \item (ii) a single exponential call, 
    \item and (iii) one multiplication of a $d\times d$ unit-norm matrix with a $d$-vector.  
\end{itemize}
Steps (i)–(ii) are $O(1)$, while (iii) is $O(d)$.
Summing over all $|E|$ edges therefore costs $|E|\!\cdot\!O(d)=\Theta(|E|d)$. This result matches the complexity of the vanilla GCN aggregation $A f$, where $A$ is the weighted adjacency matrix. Since the weight $w_{\alpha}(u,v)$ depends only on the fixed scalar $\alpha$ and pre-computed geometry, no additional learnable parameters are introduced.

Lemma~\ref{lem:gap} guarantees $\lvert\beta_{o,\xi}(u)-\beta_{o,\xi}(v)\rvert\le 1$ for every $(u,v)\in E$. For paths that deviate from a geodesic corridor, the thin-triangle
property of a $\delta$-hyperbolic graph tightens this to
$\lvert\beta_{o,\xi}(u)-\beta_{o,\xi}(v)\rvert\le\delta$.
It follows that
\[
  e^{-\alpha}\;\le\;w_{\alpha}(u,v)\;\le\;e^{-\alpha\delta},
\]
so larger curvature parameter $\delta$ pushes the weight toward the lower limit $e^{-\alpha\delta}$, yielding stronger attenuation. Subsection \ref{subsec:prove} converts this edge-wise bound into the spectral result $\lVert T_{\alpha}\rVert_{2}\le e^{-\alpha\delta}$. This matches the standard complexity analysis for a GCN layer by Kipf et al., which scales as $\Theta(|E|d)$ per layer in sparse implementation \cite{kipf2017gcn}.

\smallskip
\textbf{Multi-anchor variant.}
In many graphs, a single boundary ray is too restrictive because distinct branches may diverge along different directions. To adapt to such diversity, we select a finite set of boundary points
$\{\xi_{1},\dots,\xi_{M}\}\subset\partial G$
and assign non-negative mixing coefficients
$\alpha_{m}$ that obey $\sum_{m=1}^{M}\alpha_{m}=1$.

For each edge $(u,v)$, we replace the single-anchor weight in \eqref{eq:edgeweight} by the convex combination
$$
  w_{\boldsymbol{\alpha}}(u,v)
  \;=\;
  \sum_{m=1}^{M}\alpha_{m}\,
      \exp\!\Bigl[-\alpha_{m}\,
        \bigl|\beta_{o,\xi_{m}}(u)-\beta_{o,\xi_{m}}(v)\bigr|\Bigr].
$$
Each term quantifies the directional deviation of the edge $(u,v)$ with respect to a single boundary ray. As a result, the total weight becomes small only if divergence occurs across all selected rays.

Let $\bar\alpha=\max_{m}\alpha_{m}$. Since each exponential factor satisfies
$e^{-\bar\alpha\delta}\le\exp[-\alpha_{m}\,|\beta_{o,\xi_{m}}(u)-\beta_{o,\xi_{m}}(v)|]\le 1$,
the combined weight is bounded by
$e^{-\bar\alpha\delta}\le w_{\boldsymbol{\alpha}}(u,v)\le 1$.
Therefore, every operator-norm result proven for the single-anchor case holds \emph{verbatim} after substituting $\alpha\mapsto\bar\alpha$.

The propagation rule becomes
$$
  (T_{\boldsymbol{\alpha}}f)(v)=
  \sum_{u\in N(v)} w_{\boldsymbol{\alpha}}(u,v)\,A\,f(u),
$$
and the per-edge overhead is $M$ cached exponentials; typically $M\!\ll\!|E|$ in practice.

\begin{figure}[h]
  \centering
  \includegraphics[width=\linewidth]{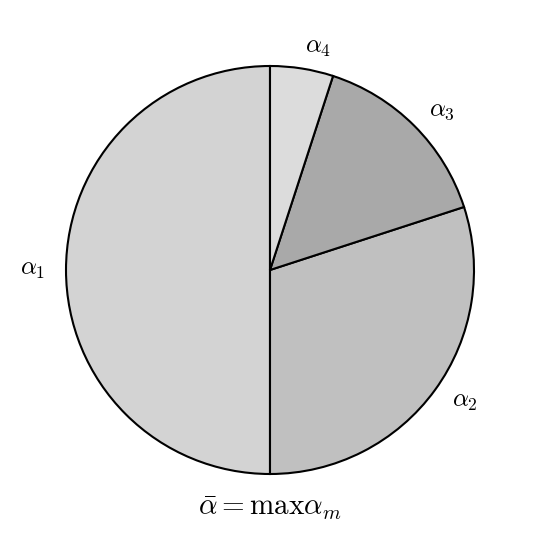}%
  \caption{Convex mixture of $M$ boundary directions.  Each slice represents a coefficient $\alpha_m$. The largest slice equals $\bar{\alpha}$, which controls the worst-case spectral bound.}
  \label{fig:multi_anchor}
\end{figure}

Figure~\ref{fig:multi_anchor} visualises the convex weights
$\{\alpha_m\}_{m=1}^M$ that modulate the mixture. The largest slice corresponds to $\bar{\alpha}$ and therefore dictates
the worst-case contraction factor in Corollary 1. When all $\alpha_m$ are similar, the pie is nearly uniform. It means
that no single boundary direction dominates the stability bound.

\subsection{Spectral guarantees}\label{subsec:prove}

We now formalize the theoretical guarantees for the propagation scheme introduced in the prior section.

\subsubsection{Primary theorem}

\begin{thm}[Boundary-Weighted Contraction Theorem]\label{thm:contraction}

Let $T_\alpha$ be the operator with row-stochastic weights $\widetilde{w}_\alpha(u,v)$. This follows that
\[
\|T_\alpha\|_2\;\le\;1
\]
for every $\alpha>0$.  

If the normalised weights are used instead and
$\alpha>\log\Delta$, one has the strict bound
\(
\|T_\alpha\|_2\le e^{-(\alpha-\log\Delta)}.
\)
\end{thm}

In order to prove this theorem, we prepare the following lemma.

\medskip
\subsubsection{Preparatory lemma}

\begin{lem}[Edge-wise Busemann Lipschitz bound]\label{lem:gap}
For every edge $(u,v)\in E$,
\[
\bigl|\beta_{o,\xi}(u)-\beta_{o,\xi}(v)\bigr| \;\le\; 1,
\]
and equality holds iff $(u,v)$ lies on some geodesic ray asymptotic to $\xi$.
\end{lem}

\begin{figure}[h]
  \centering
  \includegraphics[width=\linewidth]{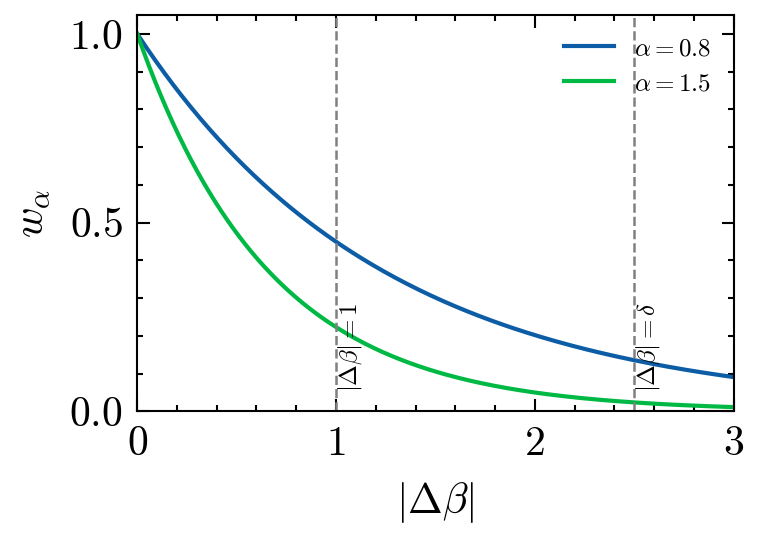}
  \caption{Edge-weight attenuation is shown. The boundary factor $w_{\alpha}$ decays exponentially with the Busemann gap $|\Delta\beta|$.  Dashed lines indicate the generic bounds $|\Delta\beta|=1$ and $|\Delta\beta|=\delta$ used in the proof of Lemma~\ref{lem:gap}.}
  \label{fig:attenuation}
\end{figure}

\noindent
Lemma~\ref{lem:gap} implies that the boundary weight $w_{\alpha}$ is sandwiched between $e^{-\alpha}$ and $e^{-\alpha\delta}$. This behavior is shown in fig.~\ref{fig:attenuation}.

\begin{proof}
We start by writing $\beta(x)=\lim_{t\to\infty}\bigl(d(x,\gamma_\xi(t))-t\bigr)$.

For any fixed $t$,
\(
|d(u,\gamma_\xi(t))-d(v,\gamma_\xi(t))|\le d(u,v)=1
\)
by the triangle inequality. Taking the limit as $t\to\infty$ preserves the inequality, giving the desired bound. 

If $(u,v)$ itself sits on the ray $\gamma_\xi$, then $d(u,\gamma_\xi(t))=d(v,\gamma_\xi(t))\pm1$ for all large $t$, so the limit equals~1. Conversely, if the limit equals~1, standard arguments on gradient lines of Busemann functions in hyperbolic spaces imply that $(u,v)$ lies on some geodesic asymptotic to $\xi$ \cite{Sormani98, DalboPeigneSambusetti12}.
\end{proof}

\medskip
\subsubsection{Proof of Theorem \ref{thm:contraction}}
\begin{figure}[h]
  \centering
  \includegraphics[width=\linewidth]{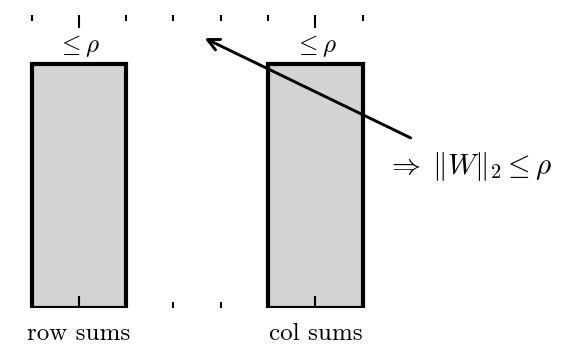}
  \vspace{-10pt}
  \caption{Row and column sums are both bounded by
           $\rho$\,\,$\Rightarrow$\. The induced $2$-norm satisfies $\|W\|_{2}\le\rho$.}
  \label{fig:contraction_sketch}
\end{figure}

\begin{proof}
Lemma~\ref{lem:gap} gives the \emph{upper} bound
\(|\beta_{o,\xi}(u)-\beta_{o,\xi}(v)|\le 1\) for every edge
$(u,v)\!\in\!E$. As a result,
\[
   w_\alpha(u,v)\;=\;e^{-\alpha|\beta(u)-\beta(v)|}\;\le\;e^{-\alpha}.
\]

Here we let \(\Delta=\max_{v\in V}\deg(v)\). For each vertex \(v\) the row sum follows
\(
 s(v)=\sum_{u\in N(v)}w_\alpha(u,v)\le\deg(v)e^{-\alpha}\le\Delta e^{-\alpha}.
\) The same estimate holds for every column sum
\(c(u)=\sum_{v\in N(u)}w_\alpha(u,v)\) because the graph is simple and undirected.  

Figure~\ref{fig:contraction_sketch} shows this estimate. In the figure, both $\|W\|_\infty$ and $\|W\|_1$ bounded by
\(\rho := \Delta e^{-\alpha}\). The left bar represents the maximum row sum $\|W\|_\infty$, the right bar the maximum column sum $\|W\|_1$, and both bars are capped at the common level $\rho$.  The diagonal arrow shows that the standard induced-norm inequality $\|W\|_2\le\sqrt{\|W\|_1\|W\|_\infty}$ merges these two one-dimensional bounds into the spectral contraction used in the next line of the proof. 

As stated previously, we define
\(
   \rho := \Delta e^{-\alpha}.
\) It follows that
\(
  \|W\|_\infty = \max_v s(v)\le\rho,\qquad
  \|W\|_1      = \max_u c(u)\le\rho.
\)  
(The row/column‐sum norms appear, e.g., in \cite[§5.2]{Meyer00} and guarantee the next step.)

The sub‐multiplicative relationship between induced norms
\begin{equation}
   \|A\|_2 \le \sqrt{\|A\|_1 \|A\|_\infty}
\end{equation}
now returns
\[
   \|W\|_2 \le \sqrt{\rho\, \rho} = \rho.
\]

If \(\alpha>\log\Delta\) then \(\rho<1\), so \(W\) is a strict
$\ell^2$-contraction.  For \(\alpha=\log\Delta\) it is non-expansive. Therefore, for any signal $x\in\mathbb R^{|V|d}$,
\[
   \|T_{\alpha}x\|_{2}
   \;=\;\|Wx\|_{2}
   \;\le\;\rho\,\|x\|_{2}
   \;=\;\Delta e^{-\alpha}\|x\|_{2}
\]

This completes the proof.
\end{proof}

\section{Conclusion}\label{sec:concl}
We have so far introduced a boundary–weighted graph filter that incorporates Gromov–boundary information through a single Busemann‐based edge weight.  For any finite $\delta$-hyperbolic graph, the filter is provably $\ell^{2}$-contractive, with an explicit curvature–controlled bound
$\|T_{\alpha}\|_{2}\!\le\!e^{-\alpha\delta}$.  The result links negative curvature, boundary geometry, and spectral stability in closed form while preserving the linear $\,\Theta(|E|d)$ complexity of a standard graph convolution and requiring no trainable parameters.

Our proof technique extends to (i)~stacked layers, which suggests the exponential decay $\|T_{\alpha}^{\,k}\|_{2}\!\le\!
e^{-k\alpha\delta}$, and (ii)~convex mixtures of multiple boundary rays, which remain contractive under the same bound with $\alpha\!\mapsto\!\bar{\alpha}$.  Future work will examine nonlinear activations and infinite graphs with bounded growth. We aim to turn the present linear theory into a full stability analysis for boundary-aware deep architectures.

\newpage


\end{document}